\documentclass[12pt]{article}
\usepackage{amsfonts,amssymb,amsbsy,amsmath,latexsym,natbib,graphicx}
\usepackage{psfrag,epsfig,epsf} 
\usepackage[margin=1in]{geometry}
\usepackage[usenames,dvipsnames,svgnames,table]{xcolor}
\usepackage{braket}

\everymath{\color{OliveGreen}}

\newcommand{\cola}[1]{\textcolor{blue}{#1}} 
\newcommand{\colc}[1]{\textcolor{red}{#1}} 

\newcommand{\Prob}[1]{\mathbb{P}\left({#1}\right)}
\newcommand{\Expect}[1]{\mathbb{E}\left[{#1}\right]}
\newcommand{\Expects}[2]{\mathbb{E}_{{#1}}\left[{#2}\right]}
\newcommand{\md}{\mbox{d}}

\newcommand{\sfX}{\mathsf{X}}

\newcommand{\setA}{\mathcal{A}}

\newcommand{\setH}{\mathcal{H}}
\newcommand{\alphabar}{\overline{\alpha}}

\newcommand{\Cov}[1]{\mbox{Cov}\left[{#1}\right]}
\newcommand{\Var}[1]{\mbox{Var}\left[{#1}\right]}
\newcommand{\Vars}[2]{\mbox{Var}_{{#1}}\left[{#2}\right]}
\newcommand{\Dir}{\mathcal{E}}
\newcommand{\mudxy}{\mu(\md x,\md y)}
\newcommand{\lambdabar}{\overline{\lambda}}
\newcommand{\cH}{\mathcal{H}}

\begin{document}
\begin{center}
\begin{Large}
\cola{\textbf{Reversible Markov chains}}
\end{Large}

\cola{Variational representations and ordering}

\textit{Chris Sherlock}

\cola{\textbf{Abstract}}

\small{This pedagogical document explains three variational representations that
  are useful when comparing the efficiencies of reversible Markov
  chains: (i) the Dirichlet form and the associated variational
  representations of the spectral gaps; (ii) a variational representation of
  the asymptotic variance of an ergodic average; and (iii) the
  conductance, and the equivalence of a non-zero conductance to a
  non-zero right spectral gap.}

\end{center}

\subsection*{\cola{Introduction}}
This document relates to variational representations of aspects of a
reversible Markov kernel, $P$, which has a limiting (hence, stationary)
distribution, $\pi$. It is a record of my own learning,
but I hope it might be useful to others. The central idea is the \colc{Dirichlet
  form}, which can be thought of as a generalisation of expected
squared jumping distance to cover all functions that are square-integrable with
respect to $\pi$. The
minimum (over all such functions with an expectation of $0$ and a variance of $1$) of the Dirichlet forms  is the right
spectral gap. 

A key quantity of interest to researchers in MCMC is the
\colc{asymptotic variance} which relates to the  
the variance of an ergodic average,
$\Var{\frac{1}{n}\sum_{i=1}^nf(X_i)}$ as $n\rightarrow \infty$, and where $X_1,\dots,X_n$ arise
from successive applications of $P$, and $f$ is some square-integrable
function. A second variational representation, of $\braket{f,(I-P)^{-1}f}$, allows
us to relate a limit of the variance of an ergodic average 
to the Dirichlet
form. Finally, a third variational quantity, the \colc{conductance} of a Markov kernel is introduced,
and is also related to the Dirichlet form and hence to the variance of
an ergodic average. 

I found these variational representations particular useful in two
pieces of research that I performed, mostly in 2016. 
\cite{SherlockThieryLee:2017} uses the first two representations,
while 
 \cite{SherlockLee:2017} uses
conductance. In both pieces of work the variational representations
allowed us to compare pairs of Markov kernels; if one Markov kernel has a
particular property, such as the variance of an ergodic average being
a particular value, and if we can relate  aspects of this Markov
kernel, such as its Dirichlet form or its conductance, to a second
Markov kernel, then we can often obtain a bound on the property of
interest for the second kernel. This is not the only use of
variational representations; e.g. in \cite{LawlerSokal:1988}
conductance is used directly to obtain bounds on the spectral gap of
several discrete-statespace Markov chains.

The most natural framework for representing the kernel 
is that of a bounded, self-adjoint operator on a Hilbert space. I was almost entirely unfamiliar
with this before I embarked on the above pieces of research and so I
will start by setting down the key aspects of this.

\subsection*{\cola{Preliminaries}}
\subsubsection*{\cola{Hilbert Space}}
Let $L^2(\pi)$ be the Hilbert space of real functions that are square
integrable with respect to some probability measure, $\pi$:
\begin{align*}
\int \pi(\md x) f(x)^2<\infty \Leftrightarrow f\in L^2(\pi),
\end{align*}
equipped with the inner product (finite through Cauchy-Schwarz) and associated norm:
\begin{align*}
\braket{f,g}=\int \pi(\md x) f(x) g(x),
~~~
||f||^2 = \braket{f,f}.
\end{align*}
Let $L_0^2(\pi)\subset L^2(\pi)$ be the Hilbert space that uses the
same inner product but includes only functions with $\Expects{\pi}{f}=0$:
\begin{align*}
\int \pi(\md x) f(x) = \braket{f,1} = 0.
\end{align*}
For these functions, $\braket{f,g}=\Cov{f,g}$ and $\braket{f,f}=\Vars{\pi}{f}$. 
\subsubsection*{\cola{Markov kernel and detailed balance}}
Let $\{X_t\}_{t=0}^\infty$ be a Markov chain on a statespace $\sfX$
with a kernel of $P(x,\md y)$ which satisfies detailed balance with
respect to $\pi$:
\begin{align*}
\pi(\md x) P(x,\md y) = \pi(\md y) P(y,\md x).
\end{align*}
For any measure, $\nu$, we define
\begin{align*}
\nu P:=\int_x \nu(\md x) P(x,\md y)
\end{align*}
Then
\begin{align*}
\pi P= \int_x \pi(\md x) P(x,\md y)
=\int_x \pi(\md y)P(y,\md x) = \pi, 
\end{align*}
so $\pi$ is stationary for $P$.

\subsubsection*{\cola{$P$ is bounded and self adjoint}}
Given a kernel (or `operator') $P$ we use the shorthand:
\begin{align*}
Pf(x)=(Pf)(x):=\int P(x,\md y) f(y).
\end{align*}
Jensen's inequality gives
\begin{align*}
[(Pf)(x)]^2=
\left[\int P(x,\md y)f(y)\right]^2 \le
\int P(x,\md y)f(y)^2
=(Pf^2)(x).
\end{align*}
Hence, because $\pi$ is stationary for $P$,
\begin{align*}
||Pf||^2 =\int \pi(\md x)[(Pf)(x)]^2 \le \int \pi(\md x) (P f^2)(x)
=
\iint 
\pi(\md x) P(x,\md y) f^2(y)
=
\int \pi(\md y) f^2(y)
=
||f||^2.
\end{align*}
Thus $||Pf||/||f|| \le 1$, and $P$ is a \colc{bounded} linear
operator.

Further, if $P$ satisfies detailed balance with respect to $\pi$
\begin{eqnarray*}
\braket{Pf,g}&=&\iint \pi(\md x) P(x,\md y) f(y) g(x)
= 
\iint \pi(\md y) P(y,\md x) f(y) g(x)\\
&=&
\iint \pi(\md x) P(x,\md y) f(x) g(y)
=\braket{f,Pg};
\end{eqnarray*}
$P$ is \colc{self-adjoint}.

\subsubsection*{\cola{The spectrum of a bounded, self--adjoint operator}}
The \colc{spectrum} of $P$
in $\setH$ is 
$\{\rho: P-\rho I~\mbox{\textcolor{black}{is not invertible in }}\setH\}$\footnote{i.e. there is
  at least one $g\in \setH$ such that there is no \textit{unique} $f\in \setH$
  with $(P-\rho I)f=g$}; the spectrum of a bounded, self-adjoint
operator is (accept it, but see below) a closed, bounded set on the real line; let the
  upper and lower bounds be $\lambda^{\max}$ and $\lambda^{\min}$. 
When $P$ is a self-adjoint Markov kernel and $\setH=L^2(\pi)$,
$\lambda^{\max}=1$ and $\lambda^{\min}\ge -1$. The
  spectral decomposition theorem for bounded, self-adjoint operators states
  that there is a finite positive measure, $a^*(\md \lambda)$ with
  support contained in the real interval
  $[\lambda^{\min},\lambda^{\max}]$ such that 
\begin{align*}
\braket{f,P^nf}=\int_{\lambda^{\min}}^{\lambda^{\max}}\lambda^n a^*(\md \lambda).
\end{align*} 
This decomposition is used twice hereafter; however it may be unfamiliar 
 so the remainder of this section provides an intuition in terms of
 the eigenfunctions and eigenvalues of $P$.

Let $P$ be a bounded, self-adjoint operator. A right eigenfunction
$e$ of $P$ is a function that satisfies 
$Pe = \lambda e$
for some scalar, $\lambda$, the corresponding eigenvalue. 
Let $e_0(x),e_1(x),\dots$ be a set of eigenfunctions of $P$, scaled so that
$||e_i||^2=\braket{e_i,e_i}=1$, and with corresponding eigenvectors
$\lambda_0,\lambda_1,\dots$. 
Since, by definition, $(P-\lambda_i I)e_i=0$, the spectrum is a
superset of the set of eigenvalues. The intuition below comes from the case
where the spectrum is precisely the set of eigenvalues, and the
eigenfunctions span $\setH$.
\begin{enumerate}
\item
Just as for the eigenvectors of a finite
self-adjoint matrix, it is possible to choose the eigenfunctions such that
$\braket{e_i,e_j}=0~(i\ne j)$; 
\item
moreover, as with self-adjoint
matrices, all of the eigenvalues are
real. 
\item
Furthermore, since $P$ is bounded, all of the
eigenvalues satisfy $|\lambda|\le 1$.
\end{enumerate}
If the eigenfunctions of an operator, $P$, span the Hilbert
space, $\setH$, for any $f\in \setH$,
\begin{align*}
f=\sum_{i=0}^\infty a_ie_i,
\end{align*}
where $a_i=\braket{f,e_i}$, and from which
\begin{align*}
\braket{f,P^nf}=\sum_{i=0}^\infty \sum_{j=0}^\infty
a_ia_j\braket{e_i,P^ne_j}
=\sum_{i=0}^\infty a_i^2\lambda_i^n.
\end{align*}
With $n=0$ we obtain $\sum_{i=0}^\infty a_i^2 = ||f||^2<\infty$.
Thus, if the eigenfunctions span the space then $a^*$ is discrete 
with mass $a_i^2$ at $\lambda_i$. 

For a Markov kernel, $P1=1$, where $1$ is the
constant function and is, thus, a right eigenfunction with an
eigenvalue of $1$.

\subsection*{\cola{Spectral gaps and the Dirichlet form}}

\subsubsection*{\cola{Spectrum of $P$ in $L_0^2(\pi)$; spectral gaps
    and geometric ergodicity}}
For any kernel $P$, we define $P^{k}f=P^{k-1}Pf$, recursively.

Any function $f(x)\in L^2(\pi)$ can be written as $a'1+f_0(x)$ where
$f_0\in L_0^2(\pi)$. Here $a'=\braket{f,1}=\Expects{\pi}{f(X)}$. Thus
\begin{align*}
P^n f = \Expects{\pi}{f(X)} + P^nf_0(x)\rightarrow \Expects{\pi}{f(X)}
\end{align*}
if $P$ is ergodic. To bound the size of the remainder term, 
consider the spectrum of $P$ restricted to functions in $L_0^2(\pi)$.
This must be confined to $[\lambda_0^{\min},\lambda_0^{\max}]$ with 
$-1 \le \lambda_0^{\min}\le \lambda_0^{\max}\le 1$, where (by definition)
\begin{align*}
\lambda_0^{\max}:=\sup_{f\in L_0^2(\pi)}\frac{\braket{f,Pf}}{\braket{f,f}}
=
\sup_{f\in L_0^2(\pi):\braket{f,f}=1}\braket{f,Pf}
~~~\mbox{\textcolor{black}{and}}~~~
\lambda_0^{\min}:=
\inf_{f\in L_0^2(\pi):\braket{f,f}=1}\braket{f,Pf}.
\end{align*}
Let $\lambdabar=\max(|\lambda_0^{\min}|,|\lambda_0^{\max}|)$ and let $a_0^*$ be the appropriate measure for $f_0$ and $P$. The (squared) size of the remainder is then
 \begin{align*}
||P^nf_0||^2=\braket{P^nf_0,P^nf_0}=\braket{f_0,P^{2n}f_0}=
\int_{\lambda_0^{\min}}^{\lambda_0^{\max}}
\lambda^{2n} a_0^*(\md \lambda)
\le 
\lambdabar^{2n}\int_{\lambda_0^{\min}}^{\lambda_0^{\max}}
a_0^*(\md \lambda)
=
\lambdabar^{2n}||f_0||^2.
\end{align*}
Thus $||P^nf_0||\rightarrow 0$ geometrically quickly provided
$\lambdabar<1$; i.e., provided the inequality is strict. $P$ is then
called \colc{geometrically ergodic}.
The \colc{right spectral gap} is $\rho^{\mathrm{right}}:=1-\lambda_0^{\max}$ and the \colc{left spectral gap}
is $\rho^{\mathrm{left}}:=1+\lambda_0^{\min}$. Both must be non-zero for geometric ergodicity.
Henceforth, for notational simplicity, we drop the subscript in $\lambda_0^{\max/\min}$.

\textbf{Aside}: when one or both of the spectral gaps is zero
(e.g. the spectrum is $[-1,1]$), then for any fixed $n$ we can always
find functions, $f$, with $||f||=1$, (albeit `fewer and fewer' as
$n\rightarrow \infty$) such that $||P^nf||>0.1$, say. 

\subsubsection*{\cola{Dirichlet form, $\Dir_P(f)$}}
The concept of a spectral gap motivates the \colc{Dirichlet form} for $P$ and $f$,
\begin{align*}
\Dir_P(f):=\braket{f,(I-P)f}=\braket{f,f}-\braket{f,Pf},
\end{align*}
since 
\begin{align*}
\rho^{\mathrm{right}}= \inf_{f\in
  L_0^2(\pi):\braket{f,f}=1}\Dir_P(f)
~~~\mbox{and}~~~
\rho^{\mathrm{left}}= 2-\sup_{f\in
  L_0^2(\pi):\braket{f,f}=1}\Dir_P(f).
\end{align*}
Directly from the definition we have for two kernels $P_1$ and $P_2$
and some $\gamma>0$:
\begin{align}
\label{eqn.gap.ordering}
\mathcal{E}_{P_1}(f)\ge \gamma \mathcal{E}_{P_2}(f)~\forall~f\in
L_0^2(\pi)\Rightarrow 
\rho^{\mathrm{right}}_{P_1}\ge \gamma \rho^{\mathrm{right}}_{P_2}.
\end{align}
An alternative expression for the Dirichlet form provides a very natural intuition:
\begin{eqnarray}
\nonumber
\Dir_P(f)&=&\braket{f,f}-\braket{f,Pf}\\
\nonumber
&=& \iint\pi(\md x)P(x,\md y) f(x)[f(x)-f(y)]\\
\nonumber
&=& \iint\pi(\md x)P(x,\md y) f(y)[f(y)-f(x)]\\
\label{eqn.Diri.alternative}
&=&
\frac{1}{2} \iint \pi(\md x)P(x,\md y)[f(y)-f(x)]^2,
\end{eqnarray}
where the penultimate line follows because $P$ satisfies detailed
balance with respect to $\pi$ and the final line arises from the average of the two
preceding lines. The Dirichlet form can, therefore, be thought of as a
generalisation of expected squared jumping distance of the $i$th
component of $x$, $\int \pi(\md x) P(x,\md y) (y_i-x_i)^2$, to consider the expected squared changes for any $f \in L^2(\pi)$.

\subsection*{\cola{Variance of an ergodic average}}
Suppose that we are interested in $\Expects{\pi}{h(X)}$ for some 
$h\in L^2(\pi)$ and we estimate it by an average of the values in the
Markov chain:
$\hat{h}_n:=\frac{1}{n}\sum_{i=1}^nP^{i-1}h(x)$. Typically  
$\Var{\hat{h}_n}\downarrow 0$ as $n\rightarrow \infty$, but scaling by
$\sqrt{n}$ should keep it $\mathcal{O}(1)$.

We are, therefore, interested in \begin{align*}
\mbox{Var}(P,h):=\lim_{n\rightarrow \infty}\Var{\sqrt{n} \hat{h}_n}.
\end{align*}
So as to just consider mixing, we assume $X_0\sim \pi$.

Without loss of generality we may assume $h\in L_0^2(\pi)$ (else just
subtract its expectation). 
Since $P$ is time-homogeneous,
\begin{eqnarray*}
\Var{\sum_{i=1}^nh(X_i)}
&=&\sum_{i=1}^n\Expect{h^2(X_i)} +
2\sum_{i=1}^{n-1}\Expect{h(X_i),h(X_{i+1})}+
2\sum_{i=1}^{n-2}\Expect{h(X_i),h(X_{i+2})}\\
&&+
\dots+ 2\Expect{h(X_1),h(X_{n})}\\
&=&
n\Expect{h^2(X)}+2(n-1)\Expect{h(X_1),h(X_2)}+
\dots+2\Expect{h(X_1),h(X_{n})}\\
&=&
n\braket{h,h}+2(n-1)\braket{h,Ph}+2(n-2)\braket{h,P^2h}+\dots+2\braket{h,P^{n-1}h}.
\end{eqnarray*}
Using $(I-P)^{-1}$ to denote $I+P+P^2+P^3+\dots$, we obtain \footnote{This is easy to see if the left and right spectral
gaps of $P$ are non-zero. When at least one gap is zero, for any $\beta<1$, we have 
$\mbox{Var}(\beta P,h)=\braket{h,h}+2(1-1/n)\beta\braket{h,Ph}+2(1-2/n)\beta^2\braket{h,P^2h}+\dots+2(1-(n-1)/n)\beta^{n-1}\braket{h,P^{n-1}h}\rightarrow\braket{h,h}+2\sum_{i=1}^\infty\beta^i\braket{h,P^ih}<\infty$. 
Letting
$\beta\uparrow 1$ then gives the result, even though $\mbox{Var}(P,h)$
may be infinite.}
\begin{eqnarray*}
\mbox{Var}(P,h)
&=&\braket{h,h}+2\sum_{i=1}^\infty\braket{h,P^ih}\\
&=&2\sum_{i=0}^\infty\braket{h,P^ih}-\braket{h,h}\\
&=&2\braket{h,(I-P)^{-1}h}-\braket{h,h},
\end{eqnarray*}
Writing 
$\braket{h,h}=\braket{h,(I-P)(I-P)^{-1}h}$ gives an equivalent 
form:
\begin{align*}
\mbox{Var}(P,h)= \braket{h,(I+P)(I-P)^{-1}h}.
\end{align*}
The spectral decomposition of $P$ in $L_0^2(\pi)$ gives
\begin{align}
\label{eqn.var.eigen.bound}
\mbox{Var}(P,h)=\int_{\lambda^{\min}}^{\lambda^{\max}}\frac{1+\lambda}{1-\lambda}~a^*(\md
\lambda)
\le
\frac{1+\lambda^{\max}}{1-\lambda^{\max}}\int_{\lambda^{\min}}^{\lambda^{\max}}a^*(\md
\lambda)
=
\frac{1+\lambda^{\max}}{1-\lambda^{\max}}\Vars{\pi}{h},
\end{align}
where the inequality follows since $(1+\lambda)/(1-\lambda)$ is an increasing function of
$\lambda$.
The supremum of the spectrum of $L_0^2(\pi)$ provides an efficiency bound over all $h \in L_0^2(\pi)$.

\subsubsection*{\cola{Variance bounding kernels}}
The bound on the variance in \eqref{eqn.var.eigen.bound} is in terms of $\lambda^{\max}$,
not $\lambda^{\min}$. Even if there exists an eigenvalue
$\lambda^{\min}=-1$, so the chain never
converges (as opposed to the more usual case where the spectrum of $P$
is $[-1,a]\cup {1}$, but there is no eigenvalue at $-1$) the asymptotic variance can be finite. 
As an extreme example, consider the 
Markov chain with a transition matrix of
\begin{align*}
\left(
\begin{array}{rr}
0&1\\1&0
\end{array}
\right).
\end{align*}
After $2n$ iterations this has been in state $1$ exactly
$1/2=\pi(\{1\})$ of the time. 

\cite{RobertsRosenthal:2008} realised that in almost all applications
of MCMC it was the variance of the ergodic average that was important,
and not (directly) convergence to the target. Specifically, it did
not matter if $\rho^{\mathrm{left}}=0$. A kernel where $\rho^{\mathrm{right}}>0$ was
termed \colc{variance bounding} and this was shown to be equivalent to 
$\mbox{Var}(P,f)<\infty$ for all $f\in L^2(\pi)$; i.e. if
[hypothetical!] simple Monte
Carlo would lead to the standard
$\sqrt{n}$ rate of convergence of the ergodic average, a $\sqrt{n}$ rate would
also be observed from the MCMC kernel.

\subsubsection*{\cola{Variational representation of
    $\braket{f,Q^{-1}f}$, and a key ordering}}

Firstly, notice that for real numbers $f,g$ and $q$,
$\sup_{g}2fg-g^2q$ is $f^2/q$, which is achieved when $g=f/q$. The
same holds for the operator $Q=I-P$ when $f$ and $g$ are functions
(see Appendix \ref{sec.proof.variation.inv} for a proof):
\begin{align*}
\braket{f,(I-P)^{-1}f} = \sup_{g\in
  L_0^2(\pi)}2\braket{f,g}-\braket{g,(I-P)g}
=\sup_{g\in L^2_0(\pi)}2\braket{f,g}-\mathcal{E}_P(g).
\end{align*}
To see this, think of a reversible matrix, $P$ and diagonalise it; the
result for such matrices is just equivalent to the scalar result
applied to each eigenvalue.

This leads to the alternative representation:
\begin{align*}
\mbox{Var}(P,h)=\sup_{g\in L_0^2(\pi)}4\braket{h,g}-2\Dir_P(g)-\braket{h,h}.
\end{align*}
Clearly, 
$\Dir_{P_1}(g)\ge \Dir_{P_2}(g)~\forall~g\in L_0^2(\pi)\Rightarrow 
\mbox{Var}(P_1,h) \le \mbox{Var}(P_2,h)$,
 giving a \textit{much} simpler proof of a result which was
 originally proved for finite-statespace Markov
chains in  \cite{Peskun:1973} and then generalised to general statespaces by \cite{Tierney:1998}. But we can
go further:
suppose that $\mathcal{E}_{P_1}(g)\ge \gamma \mathcal{E}_{P_2}(g)$ for all
$g\in L_0^2(\pi)$ and some $\gamma>0$, then, for any $g\in L_0^2(\pi)$
\begin{align*}
2\braket{h,g}-\mathcal{E}_{P_1}(g)\le
2\braket{h,g}-\gamma\mathcal{E}_{P_2}(g)
=
\frac{1}{\gamma}\left\{2\braket{h,g_*}-\mathcal{E}_{P_2}(g_*)\right\},
\end{align*}
where $g_*:=\gamma g\in L_0^2(\pi)$. So
\begin{align*}
\sup_{g\in L_0^2(\pi)} 2\braket{h,g}-\mathcal{E}_{P_1}(g)
\le \frac{1}{\gamma} \sup_{g\in L_0^2(\pi)} 2\braket{h,g}-\mathcal{E}_{P_2}(g),
\end{align*} 
and hence
\begin{align}
\label{eqn.var.ordering}
\mathcal{E}_{P_1}(f)\ge \gamma \mathcal{E}_{P_2}(f)~\forall~f\in
L_0^2(\pi)\Rightarrow 
\mbox{Var}(P_1,h)+\braket{h,h}\le \frac{1}{\gamma} \left\{\mbox{Var}(P_2,h)+\braket{h,h}\right\}.
\end{align}
This result appears as Lemma 32 in \cite{AndrieuLeeVihola:2017} \cite[see
also][]{CPS1990}.

\subsubsection*{\cola{Propose-accept-reject kernels}}
Many kernels consist of making a proposal, which is then either
accepted or rejected:
\begin{align*}
P(x,\md y):= q(x,\md y)\alpha(x,y)+(1-\alphabar(x))\delta(y-x),
\end{align*}
where
\begin{align*}
\alphabar(x):=\int q(x,\md y)\alpha(x,y)
\end{align*}
is the average acceptance probability from $x$. Then, from 
\eqref{eqn.Diri.alternative}, 
\begin{align*}
\mathcal{E}_P(f)
=
\frac{1}{2}\iint\pi(\md x) q(x,\md
y)\alpha(x,y)\{f(y)-f(x)\}^2.
\end{align*}
The right spectral gap is, therefore,
\begin{align*}
\rho^{\mathrm{right}}
=
\frac{1}{2}\inf_{f\in L_0^2(\pi):\braket{f,f}=1}\iint\pi(\md x) q(x,\md
y)\alpha(x,y)\{f(y)-f(x)\}^2.
\end{align*}
There is a similar
formula for the left spectral gap.

Results \eqref{eqn.gap.ordering} and \eqref{eqn.var.ordering} then
lead directly to the following: if two propose-accept-reject kernels, both
reversible with respect to $\pi$, satisfy
\begin{align*}
q_1(x,\md y)\alpha_1(x,y) \ge \gamma q_2(x,\md y)\alpha_2(x,y)~~~\forall~x,y
\end{align*}
then $\rho^{\mathrm{right}}_{P_1}\ge \gamma \rho^{\mathrm{right}}_{P_2}$ and
$\mbox{Var}(P_1,h)+\Vars{\pi}{h}\le \frac{1}{\gamma} \left\{\mbox{Var}(P_2,h)+\Vars{\pi}{h}\right\}$.

From \eqref{eqn.var.ordering}, if for some $\gamma>0$ 
$\Dir_{P_1}(f)\ge\gamma \Dir_{P_2}(f)~\forall~f\in L_0^2(\pi)$ 
and $P_2$ is variance bounding, then so is $P_1$. Unfortunately it is
rare that we can be sure of the ordering of the Dirichlet forms for
all $f$; indeed we might be sure that there is no fixed $\gamma$ for
which the ordering does hold. 
When we cannot simply resort to a uniform ordering of
Dirichlet forms then the elegant concept of 
{conductance} can come to our aid.

\subsection*{\cola{Conductance}}
Consider any measurable set, $\setA\subseteq \sfX$. The \colc{conductance} of $\setA$ is defined as
\begin{align*}
\kappa_P(\setA):=\frac{1}{\pi(\setA)}\iint_{x\in \setA,y\in \setA^c} \pi(\md x)P(x,\md y),
\end{align*}
where $\pi(\setA):=\int_\setA \pi(\md x)$. 
Loosely speaking, $\kappa_P(\setA)$ is the probability of moving to
$\setA^c$ conditional on the chain currently following the stationary
distribution truncated to $\setA$: $\pi(\md x)1_{x\in\setA}$. 
Our analysis of the properties of $\kappa_P(\setA)$ will be via the symmetric quantity
\begin{align*}\kappa^*_P(\setA)=\kappa^*_P(\setA^c):=
\frac{1}{\pi(\setA)\pi(\setA^c)}\int_{x\in \setA,y\in \setA^c} \pi(\md x)P(x,\md y)=\frac{\kappa_P(\setA)}{\pi(\setA^c)}.
\end{align*}
The first equality follows because the kernel is reversible with
respect to $\pi$, and this also implies that $\pi(\setA)\kappa_P(\setA)=\pi(\setA^c)\kappa_P(\setA^c)$. Since $\kappa_P(\setA^c)\le 1$, $\kappa_P(\setA)\le \pi(\setA^c)/\pi(\setA)$, which can be arbitrarily small. To define the \colc{conductance of the kernel} $P$ we therefore only consider sets $\setA$ with $\pi(\setA)\le 1/2$:
\begin{align*}
\kappa_P:=\inf_{\setA:\pi(\setA)\le 1/2}\kappa(\setA).
\end{align*}
No such restriction is required for:
\begin{align*}
\kappa^*_P:=\inf_{\setA} \kappa^*_P(\setA).
\end{align*}
Further, since 
$\pi(\setA)\le 1/2\Rightarrow 1/2\le \pi(\setA^c)\le 1$, we have
$\kappa_P\le\kappa^*_P\le 2\kappa_P$.

Setting $f_\setA(x)=[1_{x\in \setA}-\pi(A)]/\sqrt{\pi(A)\pi(A^c)}$ (so that
$f_\setA\in L_0^2(\pi)$ and $\braket{f_\setA,f_\setA}=1$) in the expression for the Dirichlet form \eqref{eqn.Diri.alternative} gives
\begin{align*}
\Dir_P(f_\setA)=\frac{1}{2\pi(\setA)\pi(\setA^c)}\left\{\int_{x\in\setA,y\in\setA^c}\pi(\md x)P(x,\md y)+\int_{x\in\setA^c,y\in\setA}\pi(\md x)P(x,\md y)\right\}
=\kappa^*_P(\setA)=\frac{\kappa_P(\setA)}{\pi(\setA^c)}.
\end{align*}
So
\begin{align}
\label{eqn.condA}
\rho^{\mathrm{right}}
=
\inf_{f\in L_0^2(\pi):\braket{f,f}=1}\Dir_P(f)
\le
\inf_{\setA}\Dir_P(f_\setA)
=\kappa^*_P\le 2 \kappa_P.
\end{align}
Hence, if $P$ has a right spectral gap (so it is variance bounding) then its conductance is non-zero. Amazingly, the converse is
also true: if the conductance of $P$ is non-zero then $P$ has a right
spectral gap:
\begin{align}
\label{eqn.condB}
\rho^{\mathrm{right}}\ge\frac{\kappa^{*2}_P}{2}\ge \frac{\kappa_P^2}{2}.
\end{align}
Thus, non-zero conductance is equivalent to a non-zero right spectral
gap is equivalent to variance bounding. Equations \eqref{eqn.condA}
and \eqref{eqn.condB} together are sometimes called \colc{Cheeger
  bounds}.

A proof \eqref{eqn.condB} for finite Markov chains is given in
\cite{DiaconisStroock1991}. A  more
accessible and, as far as I can see, more general, proof of the looser inequality, that  $\rho^{\mathrm{right}}\ge
{\kappa_P^2}/{8}$ is given 
in \cite{LawlerSokal:1988}. This has a
`standard bit' which 
uses the Cauchy-Schwarz inequality to obtain an inequality for $\rho^{\mathrm{right}}$, a `beautiful bit' which relates
the expression from the `standard bit' to conductance, and then an
`ugly bit', which proves that the final expression is always greater
than $0$ if $\kappa_P>0$. I will follow \cite{LawlerSokal:1988} for
the first two parts and then provide a neater solution to the final
part.
 
Denote the symmetric measure $\pi(\md x)P(x,\md y)$ by $\mudxy$. We
also set $g(x)=f(x)+c$ for some (currently) arbitrary real constant,
$c$. 

\textbf{\cola{The standard bit}}. 
Then Cauchy-Schwarz (twice) and the symmetry of $\mu$ gives:
\begin{eqnarray*}
\Expects{\mu}{|g(X)^2-g(Y)^2|}^2
&=&
\Expects{\mu}{|g(X)-g(Y)|~|g(X)+g(Y)|}^2\\
&\le&
\Expects{\mu}{\{g(X)-g(Y)\}^2}\Expects{\mu}{\{g(X)+g(Y)\}^2}\\
&\le&
\Expects{\mu}{\{f(X)-f(Y)\}^2}2\Expects{\mu}{g(X)^2+g(Y)^2}\\
&=&
4\Expects{\mu}{\{f(X)-f(Y)\}^2}\Expects{\pi}{g(X)^2}.
\end{eqnarray*}
So
\begin{align*}
\Dir_P(f) \ge \frac{1}{8\Expects{\pi}{g(X)^2}}~\Expects{\mu}{|g(X)^2-g(Y)^2|}^2.
\end{align*}
\textbf{\cola{The beautiful bit}}. 
We now relate the numerator of the above expression to the conductance. The proof is symmetrical, the first half manipulates the denominator so that conductance may be used, with the second half reversing the route of the first.

Set $\setA_t:=\{x:g(x)^2\le t\}$. Then by the symmetry of $\mu$,
\begin{eqnarray*}
\int_{\sfX\times \sfX}\mudxy |g(y)^2-g(x)^2|
&=&
2
\int_{\sfX\times \sfX}\mudxy 1_{g(x)^2<g(y)^2}\{g(y)^2-g(x)^2\}\\
&=&
2
\int_{\sfX\times \sfX}\mudxy \int_{t=0}^\infty\md t 1_{g(x)^2\le t< g(y)^2}\\
&=&
2
\int_{t=0}^\infty\md t \int_{\sfX\times \sfX}\mudxy 1_{g(x)^2\le t< g(y)^2}\\
&=&
2\int_{t=0}^\infty\md t \int_{\setA_t \times \setA_T^c}\mudxy\\
&\ge&
2\kappa^*_P \int_{t=0}^\infty\md t \int_{\setA_t\times \setA_t^c}\pi(\md x) \pi(\md y)\\
&=&
\kappa^*_P\int_{\sfX\times \sfX}\pi(\md x)\pi(\md y) |g(y)^2-g(x)^2|.
\end{eqnarray*}
So
\begin{align*}
\Dir_P(f) \ge \frac{\kappa^{*2}_P\Expects{\pi\times \pi}{|g(X)^2-g(Y)^2|}^2}{8\Expects{\pi}{g(X)^2}}.
\end{align*}
But since this is true for all $c$, we have
\begin{align*}
\rho^{\mathrm{right}}\ge \frac{\kappa^{*2}_P}{8}\inf_{f\in L_0^2(\pi):\braket{f,f}=1} \sup_c\frac{\Expects{\pi\times \pi}{|\{f(Y)+c\}^2-\{f(X)+c\}^2|}^2}{\Expects{\pi}{\{f(X)+c\}^2}}.
\end{align*}
\textbf{\cola{A less ugly bit}}. 
Since ${\Expects{\pi}{\{f(X)+c\}^2}}=1+c^2$, we need to show that
\begin{align*}
\inf_{f\in L_0^2(\pi):\braket{f,f}=1} \sup_c\frac{\Expects{\pi\times
    \pi}{|\{f(Y)+c\}^2-\{f(X)+c\}^2|}}{\sqrt{1+c^2}}
\ge 1.
\end{align*}
Let $A=f(X)$ and $B=f(Y)$ be independent and identically distributed
with an expectation of $0$ and a variance of $1$. We are interested in
\begin{align}
\label{eqn.AcBc}
\frac{\Expect{|\{A+c\}^2-\{B+c\}^2|}}{\sqrt{1+c^2}}.
\end{align}
Below, I will show
that
\begin{align}
\label{eqn.nicer.ineq.thanLS}
\Expect{|A^2-B^2|}+4\Expect{|A-B|}^2\ge 2.
\end{align}
Then, as in \cite{LawlerSokal:1988}, consider letting $c\rightarrow \infty$ or 
setting $c=0$. With the former:
\begin{align*}
\lim_{c\rightarrow \infty}\frac{\Expect{|\{A+c\}^2-\{B+c\}^2|}}{\sqrt{1+c^2}} = 2\Expect{|A-B|},
\end{align*}
so if $\Expect{|A-B|}>1/2$ we are done. If not, setting $c=0$
in \eqref{eqn.AcBc} and using \eqref{eqn.nicer.ineq.thanLS} leaves us:
\begin{align*}
\Expect{|A^2-B^2|} \ge 1.
\end{align*}
To prove \eqref{eqn.nicer.ineq.thanLS}:
\begin{eqnarray*}
\Expect{|A^2-B^2|}
&=&
\int_0^\infty\Prob{|A^2-B^2|>t}\md t
=
\int_0^\infty\Prob{A^2>t+B^2}+\Prob{B^2>t+A^2}\md t\\
&=&
2\int_0^\infty\Prob{B^2>t+A^2}\md t
=
2\Expects{A}{\int_0^\infty\Prob{B^2>t+A^2}\md t|A}.
\end{eqnarray*}
But
\begin{eqnarray*}
\int_0^\infty\Prob{B^2>t+a^2}\md t
&=&
\int_{a^2}^\infty\Prob{B^2>v}\md v
=
\int_0^\infty\Prob{B^2>v}\md v - \int_0^{a^2}\Prob{B^2>v} \md v\\
&=&
1-\int_0^{|a|} 2u\Prob{|B|>u} \md u\\
&\ge&
1-2|a|\int_0^\infty \Prob{|B|>u}
=
1-2|a|\Expect{|B|}.
\end{eqnarray*}
Combining the two end results gives $\Expect{|A^2-B^2|}\ge
2-4\Expect{|A|}\Expect{|B|}$; i.e.
\begin{align*}
\Expect{|A^2-B^2|}+4\Expect{|A|}^2\ge 2.
\end{align*}
However Jensen's inequality provides: 
$\Expect{|A-B|}=\Expect{\Expect{|A-B|}|A}\ge
\Expect{|A-\Expect{B}|}=\Expect{|A|}$,
and
\eqref{eqn.nicer.ineq.thanLS} follows.

\section*{Acknowledgements}
I am grateful to Dr. Daniel Elton for providing the proof in Appendix
\ref{sec.proof.variation.inv} and Mr. Sam Holdstock and Dr. Dootika Vats for spotting errors in earlier versions of this document.
\appendix

\section{Variational representation of $\braket{f,Q^{-1}f}$}
\label{sec.proof.variation.inv}
Let $\cH$ be a Hilbert space and let $Q$ be a positive operator on
$\cH$ (i.e., an operator that has a square root). Then for $f\in\cH$,
\begin{align*}
\braket{f,Q^{-1}f}=\sup_{g\in \cH}2\mbox{Re}(\braket{f,g})-\braket{g,Qg}.
\end{align*}
(Our Hilbert space is real, so we do not need the Re() function.)

\textbf{Proof} 
Let $\phi=Q^{-1}f$, so $f=Q\phi$. The left-hand side
is then
\begin{align*}
  \braket{Q\phi,\phi}&=\braket{Q^{1/2}\phi,Q^{1/2}\phi}=||Q^{1/2}\phi||^2\\
  &\ge||Q^{1/2}\phi||^2 - ||Q^{1/2}(\phi-g)||^2\\
  &=2\mbox{Re}(\braket{Q^{1/2}\phi,Q^{1/2}g})-||Q^{1/2}g||^2\\
  &=2\mbox{Re}(\braket{Q\phi,g})-\braket{g,Qg}\\
  &=2\mbox{Re}(\braket{f,g})-\braket{g,Qg}.
\end{align*}
The construction of the inequality shows that the supremum is achieved
at $g=Q^{-1}f$.

\bibliographystyle{apalike}

\bibliography{mcmc_theory.bib}

\begin{thebibliography}{}

\bibitem[{Andrieu} et~al., 2016]{AndrieuLeeVihola:2017}
{Andrieu}, C., {Lee}, A., and {Vihola}, M. (2016).
\newblock Uniform ergodicity of the iterated conditional {SMC} and geometric
  ergodicity of particle {G}ibbs samplers.
\newblock {\em Bernoulli}.
\newblock to appear.

\bibitem[Caracciolo et~al., 1990]{CPS1990}
Caracciolo, S., Pelissetto, A., and Sokal, A.~D. (1990).
\newblock Nonlocal monte carlo algorithm for self-avoiding walks with fixed
  endpoints.
\newblock {\em Journal of Statistical Physics}, 60(1):1--53.

\bibitem[Diaconis and Stroock, 1991]{DiaconisStroock1991}
Diaconis, P. and Stroock, D. (1991).
\newblock Geometric bounds for eigenvalues of {M}arkov chains.
\newblock {\em Ann. Appl. Probab.}, 1(1):36--61.

\bibitem[Lawler and Sokal, 1988]{LawlerSokal:1988}
Lawler, G.~F. and Sokal, A.~D. (1988).
\newblock Bounds on the {$L^2$} spectrum for {M}arkov chains and {M}arkov
  processes: a generalization of {C}heeger's inequality.
\newblock {\em Trans. Amer. Math. Soc.}, 309(2):557--580.

\bibitem[Peskun, 1973]{Peskun:1973}
Peskun, P.~H. (1973).
\newblock Optimum {M}onte-{C}arlo sampling using {M}arkov chains.
\newblock {\em Biometrika}, 60:607--612.

\bibitem[Roberts and Rosenthal, 2008]{RobertsRosenthal:2008}
Roberts, G.~O. and Rosenthal, J.~S. (2008).
\newblock Variance bounding {M}arkov chains.
\newblock {\em Ann. Appl. Probab.}, 18(3):1201--1214.

\bibitem[{Sherlock} and {Lee}, 2022]{SherlockLee:2017}
{Sherlock}, C. and {Lee}, A. (2022).
\newblock {Variance bounding of delayed-acceptance kernels}.
\newblock {\em Methodology and Computing in Applied Probability},
  24(3):2237--2260.

\bibitem[Sherlock et~al., 2017]{SherlockThieryLee:2017}
Sherlock, C., Thiery, A.~H., and Lee, A. (2017).
\newblock Pseudo-marginal metropolis–hastings sampling using averages of
  unbiased estimators.
\newblock {\em Biometrika}, 104(3):727--734.

\bibitem[Tierney, 1998]{Tierney:1998}
Tierney, L. (1998).
\newblock A note on {M}etropolis--{H}astings kernels for general state spaces.
\newblock {\em Ann. Appl. Probab.}, 8(1):1--9.

\end{thebibliography}

\end{document}